\documentclass[a4paper,11pt,onecolumn,twoside,UTF8]{article}

\usepackage{amsthm}
\usepackage{amssymb}
\usepackage{amsmath,bm}
\usepackage[]{hyperref}

\title{Generic stability of linear algebraic groups over $\C[[t]]$}
\author{Chen Ling\\Fudan University \and Ningyuan Yao\\Fudan University  }

\newtheorem{theorem}{Theorem}

\newtheorem*{Conj}{Conjecture}

\newtheorem{lemma}{Lemma}
\newtheorem{corollary}[lemma]{Corollary}
\newtheorem{fact}[lemma]{Fact}
\newtheorem{remark}[lemma]{Remark}

\theoremstyle{definition}
\newtheorem{definition}{Definition}

\newcommand{\tp}{\mathrm{tp}}
\newcommand{\Th}{\mathrm{Th}}

\newcommand{\acl}{\mathrm{acl}}

\newcommand{\Z}{\mathbb{Z}}
\newcommand{\M}{\mathbb{M}}

\newcommand{\N}{\mathbb{N}}
\newcommand{\K}{\mathbb{K}}

\newcommand{\U}{\mathbb{U}}
\newcommand{\C}{\mathbb{C}}

\newcommand{\GG}{\mathfrak{G}}

\newcommand{\cO}{{\mathcal{O}}}

\newcommand{\SL}{\mathrm{SL}}
\newcommand{\GL}{\mathrm{GL}}
\newcommand{\SO}{\mathrm{SO}}
\newcommand{\alg}{\mathrm{alg}}

\newcommand{\mm}{\mathfrak{m}}
\newcommand{\kk}{\bm{k}}    
\newcommand{\res}{\mathrm{res}}
\newcommand{\trans}{\mathrm{trans}}

\newcommand{\sq}{\subseteq}

\begin{document}
\maketitle

\begin{abstract}
Let $K$ be a henselian valued field with $\cO_K$ its valuation ring, $\Gamma$ its value group, and $\kk$ its residue field. We study the definable subsets of $\cO_K$ and algebraic groups definable over $\cO_K$ in the case where $\kk$ is algebraically closed and $\Gamma$ is a $\Z$-group. 

We first describe the definable subsets of $\cO_K$, showing that every definable subset of $\cO_K$ is either res-finite or res-cofinite (see Definition \ref{def-res-finite-cofinite}). Applying this result, we show that $\GL(n,\cO_K)$ (the invertible $n$ by $n$ matrices over $\cO_K$) are generically stable for each $n$, generalizing Y. Halevi's result, where $K$ is an algebraically closed valued field \cite{Y.Halevi}.
\end{abstract}


\section{Introduction}
The notion of generically stable types were introduced  by Hrushovski and Pillay to describe the ``stable-like'' behavior in NIP environment \cite{HP-NIP-inv-measure}. Briefly, a type over a monster model, called a global type, is generically stable if it is finitely satisfiable in and definable over every small submodel. The theory of algebraically closed valued fields, denoted by ACVF, is considered a typical ``stable-like'' NIP theory, since the non-trivial generically stable types exist. As a contrast, for the purely unstable NIP theories, say  $p$-adically closed fields ($p$CF), there is no non-trivial generically stable type. 

In this paper, we study the structures which are elementarily equivalent to $\C((t))$, the field of form Laurent series over the complex numbers.  It is natural to consider $\Th(\C((t)))$ as a mixture of the ``stable-like'' theory ACVF and the ``purely unstable'' theory $p$CF. To see this, let $K\models \Th(\C((t)))$, $K_1\models $ACVF and $K_2\models p$CF, with residue fields $\kk$, $\kk_1$ and $\kk_2$; value groups $\Gamma$, $\Gamma_{1}$ and $\Gamma_{2}$, respectively. Then $\kk\equiv \kk_{1}$ and  $\Gamma\equiv \Gamma_{2}$. On the other side, according to Ax-Kochen-Ershov\cite{Ax-Kochen, Ershov}, the theory of a henselian valued field is completely determined by the theories of its residue field and its value group if its residue field has characteristic $0$.
Let $K\equiv \C((t))$. We study the definable subsets of $\cO_K$, our first result is:
\begin{theorem}
Let $\res: K\to\kk$ be the residue map. Suppose that $K\equiv \C((t))$, $X\sq \cO_K$ a definable set. Then there is finite $Z\sq \kk $ such that either $X\sq \res^{-1}(Z)$ or $\res^{-1}(\kk \backslash Z)\sq X$.
\end{theorem}
Recall from \cite{HP-NIP-inv-measure} that a group $G$ defined in a monster model $\M$ is generically stable if there is a generically stable type for $G$, where a global type $p\in S_G(\M)$ is a generically stable type for $G$ if it is generically stable for each $g\in G$.   Y. Halevi showed that $\GL(n,\cO_K)$, the invertible $n$ by $n$ matrices over $\cO_K$, is generically stable if $K$ is an algebraically closed valued field.  In this paper, We generalize this result to the case where $K\equiv \C((t))$. Let $\K$ be a monster model of $\Th(\C((t)))$, $\cO_\K$ the valuation ring of $\K$, $\U$ be the units of $\cO_\K$, and $\GL(n,\cO_\K)$ the invertible $n$ by $n$ matrices over $\cO_\K$. Then we prove the following theorem:
\begin{theorem}
Let $G$ be $\cO_\K$, or $\U$, or $\GL(n,\cO_\K)$, then 
\begin{description}
   \item   [(i)] $G$ has a (unique) generically stable type $p$.
   \item   [(ii)] $p$ is $G$-invariant. As a consequence, $G=G^{00}$.
   \item  [(iii)] $p$ is dominated by $\res(p)$ via the residue map, where $\res(p)$ is the image of $p$ under the residue map.

\end{description}
\end{theorem}

The paper is organized as follows. For the rest of this section, we recall some basic facts around valued fields and $\mathbb{C}((t))$ and introduce notations we use. In section \ref{sec-def-val-ring}, we prove that every definable subset of the valuation ring $\cO_K$ is either res-finite or res-cofinite. In section \ref{sec-Gen-O_k-U}, we give a generically stable type $p_{trans,\K}$ of $\cO_{\K}$ and $\U$. From this onwards, in section \ref{sec-Gen-GLn}, we show that for every $n$, $\GL(n,\cO_\K)$ has a unique generically stable type that is dominated by its image under the residue map.


\subsection{Preliminaries} 
 
Let  $K$ be a field, $(\Gamma,>)$ an ordered abelian group,  and $v:K\to \Gamma\cup\{\infty\}$, where $\infty>\Gamma$. We say $(K,v)$ or $K$ is a valued field, if  $v:K\to \Gamma\cup\{\infty\}$ is an onto map satisfying the followings: 
    \begin{itemize}
        \item $v(x)=\infty$ iff $x=0$.
        \item $v(x+y)\geq \min(v(x),v(y))$.
        \item $v(xy)=v(x)+v(y)$.
    \end{itemize}
    We call $v$ its valuation map and $\Gamma$ its value group. We denote $\cO_K:=\{x\in K |\  v(x)\geq 0\}$ the valuation ring of $K$, $\mm_K:=\{x\in K |\  v(x)> 0\}$
    the unique maximal ideal of $\cO_K$, and ${\kk}=\cO_K/\mm_K$ the residue field of $K$. The residue map is the natural projection ${\res}:\cO_K\to {\kk}$. For each $x\in K$ we call  $v(x)$ the valuation of $x$.   For each $x\in \cO_K$ we call $\res(x)$ the residue of $x$. If $f(x)=a_nx^n+\cdots+a_1x+a_0$ is a polynomial with each $a_i\in \cO_K$, then by $\res(f)$ we mean the polynomial 
    \[
    \res(a_n)x^n+\cdots+\res(a_1)x+\res(a_0)\in {\kk}[x].
    \]
    We call $\cO_K$ henselian, if for any $f\in \cO_K[x]$ and $\alpha\in \cO_K$ satisfying 
    $$f(\alpha)\in\mathfrak{m}_K, f'(\alpha)\notin\mathfrak{m}_K$$ 
    there exists $a\in \cO_K$ such that $f(a)=0$ and $a\equiv \alpha \; mod \; \mathfrak{m}_K$, where $f'(x)$ is the derivative of $f(x)$. We say $K$ is a henselian valued field if $\cO_K$ is henselian.

    For any $a\in \cO_K$ and $f(x)\in \cO_K[x]$, $\res (f(a))=\res(f)(\res(a))$.
    For any subset $X$ of $\cO_K$, by $\res(X)$ we mean the set $\{\res(a)|\ a\in X\}$. It is easy to see that $\res(A\cup B)=\res(A)\cup \res(B)$ and $\res(A\cap B)\sq \res(A)\cap \res(B)$ for any $A,B\sq \cO_K$.

Let $K^*=K\backslash\{0\}$ be the multiplicative group. Then the set of $n$-th powers $P_n(K^*)=\{a^n|\ a\in K^n\}$ is a subgroup of $K^*$, which is definable in the language of rings.

\begin{fact}\label{fact-nth-power}
    Let $K$ be a henselian valued field with ${\kk}$ the residue field. Suppose that ${\kk}$ is algebraically closed and $\mathrm{char}({\kk})=0$. For any $c\in K$, if $v(c)=0$, then $c\in P_n(K^*)$ for any $n\in\mathbb{N}^{>0}$.
\end{fact}
\begin{proof}
    We consider the polynomial $f(x)=x^n-c$ where $c\in K$ and $v(c)=0$, denote $\mm_K$ as the maximal ideal on the valuation ring of $K$. As ${\kk}$ is algebraically closed,  and $\res(c)$ is nonzero, $x^n-\res(c)$ has a nonzero root in ${\kk}$. So there is $b\in K$ such that $v(b)=0$ and $\res(f(b))=(\res(b))^n-\res(c)=0$.

    Thus $f(b)=b^n-c\in\mm_K$. Since $\text{char}({\kk})=0$, the derivative $f'(b)=nb^{n-1}\notin\mm_K$. Now $K$ is henselian, so $f(x)$ has a root in $b+\mm_K$, and thus $c\in P_n(K^*)$. This completes the proof.
\end{proof}

\begin{corollary}\label{cor-nth-power}
     Let $K$ be a henselian valued field, ${\kk}$ the residue field of $K$, and $\Gamma$ the value group of $K$. If ${\kk}$ is algebraically closed and $char({\kk})=0$. Then for any $b\in K^*$,  $b\in P_n(K^*)$ iff $v(b)\in n\Gamma=\{n\gamma|\ \gamma\in\Gamma\}$. In particular, $K^*/P_n(K^*)\cong \Gamma/n\Gamma$.
\end{corollary} 
\begin{proof}
    Suppose that $b\in P_n(K^*)$, then there is $a\in K^*$ such that $b=a^n$. So $v(b)=v(a^n)=nv(a)\in n\Gamma$.

    Conversely, suppose that $b\in K^*$ such that $v(b)=nv(a)$ for some $a\in K^*$. Then $v(ba^{-n})=v(b)-nv(a)=0$. By Fact  \ref{fact-nth-power}, we see that $ba^{-n}\in P_n(K^*)$. Since $a^{-n}\in P_n(K^*)$, we have $b\in P_n(K^*)$ as required. The ``in particular'' part follows directly.
\end{proof}

\begin{fact}\label{fact-lift}
   Let $K$ be a henselian valued field, with $\cO_K$ the valuation ring, ${\kk}$ the residue field of $K$. If $\mathrm{char}({\kk})=0$, then ${\kk}$ can be lifted. Namely, there is a subfield $E$ of $\cO_K$ such that $\res: E\to {\kk}$ is an isomorphism. So we can consider ${\kk}$ as a subfield of $\cO_K$.
\end{fact}


\subsection{Model theory of $\mathbb{C}((t))$}
Let  $\mathbb{C}((t)$ be the field of formal Laurent series over the field $\C$ of complex numbers.  Elements of $\mathbb{C}((t)$ are of the form $\sum_{i=n}^{\infty}a_it^i$,  where $t$ is a variable,   $n\in \mathbb{Z}$, each $a_i\in\mathbb{C}$.  $\mathbb{C}((t))$ is a henselian valued field with the valuation map: 
 \[
 v:\mathbb{C}((t))\to \Z\cup\{\infty\}, \sum_{i=n}^{\infty}a_it^i\mapsto  \min\{i |\  i\geq n \ \text{and}\  a_i\neq 0\}.
 \]
 The valuation ring  and residue field of $\C((t))$ are $\mathbb{C}[[t]]=\{\sum_{i=n}^{\infty}a_it^i|\ n\in\N\}$ and $\C$ respectively (see \cite{Dr1} for details).  

Let $L_{ring}=\{0,1,+,\times\}$ be the language of rings and 
\[
L_{VR}=L_{ring}\cup\{ \; |\; , N\}\cup\{P_n |\  n\in\mathbb{N}^{>0}\}
\]
an expansion of $L_{ring}$ for the valuation rings, where the new predicates $x|y$ is interpreted as $v(x)\leq v(y)$, $N(x)$ is interpreted as $v(x)=1$,   and $P_n(x)$ is interpreted as the set of $n$-th powers, respectively. The theory of $\C((t))$ has NIP and quantifier elimination in the language $L_{VR}$ (see \cite{De1}).

Let $K\equiv \C((t))$ and $\psi(x)$ be a one-variable atomic $L_{VR}$-formula with parameters from $K$, then $\psi$ is of one of the following four types:
\begin{description}
\item [Type (i)] $f(x)=0$, where $f\in K[x]$ is a polynomial over $K$.
\item  [Type (ii)] $v(f(x))\leq v(g(x))$ (or $f(x)|g(x)$), where $f,g\in K[x]$ are polynomials over $K$.   
\item [Type  (iii)] $P_n(f(x))$, where $f\in K[x]$ is a polynomial over $K$.
\item [Type  (iv)] $N(f(x))$, where $f\in K[x]$ is a polynomial over $K$  
\end{description}

\begin{remark}\label{rmk-atomic formulas}
    Since  $\mathbb{C}((t)) \models  N(x)\leftrightarrow (v(x)=v(t))$, we have $K\models \exists y(N(x)\leftrightarrow (v(x)=v(y)))$, which means that there is $a\in K$ such that $K\models N(x)\leftrightarrow (v(x)=v(a))$. So any atomic one-variable $L_{VR}(K)$-formula is a Boolean combination of the formulas of \textbf{Type (i)-(iii)}. 
\end{remark}

\begin{remark}\label{res-is-defianble}
    Let $K$ be a valued field with $\cO_K$ its value group, then the relation ``$\res (x)=\res(y)$'' on $\cO_K$ is definable in the language $L_{VR}$, in fact, $\res (x)=\res(y)$ is defined by the formula 
    \[
    v(x)\geq 0\wedge v(y)\geq 0\wedge v(x-y)>0.
    \]
\end{remark}

\begin{lemma}\label{lemma-f(res) is definable}
    Let $K$ be a valued field with $\cO_K$ its value group and $\kk$ its residue field. Let $f(x)\in {\kk}[x]$, then the set
    \[
    \{a\in \cO_K|\ f(\res(a))=0\}
    \]
    is definable in the language $L_{VR}$ with parameters from $\cO_K$.
\end{lemma}
\begin{proof}
 Let $b_1,\cdots,b_m$ be the roots of $f$ in $\kk$. Let $a_1,\cdots,a_m\in\cO_K$ such that $\res(a_i)=b_i$. then   
\[
\{a\in \cO_K|\ f(\res(a))=0\}=\bigcup_{i=1}^m\{a\in \cO_K|\ \res(a)=\res(a_i)\}.
\]
Since  $\{a\in \cO_K|\ \res(a)=\res(a_i)\}$ is definable for each $i\leq m$, we see that 
$\{a\in \cO_K|\ f(\res(a))=0\}$ is definable.
\end{proof}

According to Lemma \ref{lemma-f(res) is definable}, for each $f\in{\kk}[x]$, it is reasonable to consider $f(\res(x))=0$ as a  formula in the language $L_{VR}$ with parameters from $\cO_K$.

Let $K\equiv \C((t))$, then by Corollary \ref{cor-nth-power}, any   $a\in K^*$ is an $n$-th power iff $v(a)\in n\Gamma$ and $|K^*/P_n(K^*)|=|\Gamma/n\Gamma|=|\Z/n\Z|=n$.


\subsection{Notations}
 $T$ will denote   the complete theory of $\mathbb C((t))$ in the language of $L_{VR}$. Let us ﬁx a monster model $\K$ of $T$, namely $\K$ is saturated and homogeneous of cardinality $\kappa$, where $\kappa$ is very big, say strongly inaccessible.  We call an object  ``small'' or ``bounded'' if it is of cardinality $<\kappa$.  

For $A$ a subset of $\K$, an $L_{VR}(A)$-formula is a formula with parameters from $A$. If $\phi(\bar x)$ is an $L_{VR}(\K)$-formula and $A\sq \K$, then $\phi(A)$ is the collection of the realizations of $\phi(\bar x)$ from $A$, namely, $\phi(A)=\{\bar a\in A^{|\bar x|}|\ \K\models \phi(\bar a)\}$.

Let $K\prec \K$. When we speak of a set $X$ definable in $K$ or $\K$, we mean the obvious thing. In this case, we use $X(A)$ to denote the set $\varphi(A)$. When we speak of a set $X$ definable over $K$, we mean a set definable in $\K$ defined by an $L_{VR}(K)$-formula $\varphi(\bar x)$.  
If $X\sq K^n$ is definable in $K$, we sometimes use $X(x)$ to denote the formula which deﬁnes $X$ and $X(\K)$ the subset of $\K^n$ defined by the formula $X(x)$. 
In general, when we speak of a definable object (set, or group) we mean a definable object in $\K$. Let $X$ be a set definable in $K$, then  $S_X(K)$ denotes the space of complete types over $K$  concentrating on $X$.  


Let $\GL(n,\K)\sq \K^{n\times n}$  be the general linear algebraic group over $\K$, consisting of all $n\times n$ invertible matrices over $\K$. $\GL(n,\K)$ is defined by the formula $\det(\bar x)\neq 0$, hence is definable.  For convenience, we denote the formula ``$\det(\bar x)\neq 0$'' by $\GL(n,\bar x)$. 
Clearly, for any field $F$, $\GL(n, F)$ is a  group definable in $F$, defined in the language of rings.  

Let $\cO_{\K}=\{a\in \K|\ v(a)\geq 0\}$ be the valuation ring of $\K$. By speaking of a set $X$ defined on $\cO_{\K}$, we mean $X$ is definable and $X\sq \cO_{\K}^n$. By $\GL(n,\cO_{\K})$, we mean the set $\{\bar g\in \GL(n,\K)|\ v(\det(g))=0\}$. It is easy to see that $\GL(n,\cO_{\K})$ is a subgroup of $\GL(n,\K)$. Note that $\GL(n,\cO_{\K})\neq \GL(n,\K)\cap \cO_{\K}^{n\times n}$ since the latter is NOT a group.  Let $\U=\{a\in \K|\ v(a)=0\}$ be the set of units in $\cO_{\K}$, it is easy to see that $\U=\GL(1,\cO_{\K})$. It is natural to consider $\cO_{\K}$ and $\U$ as counterparts of the additive group and multiplicative group in $\cO_{\K}$. 

We use ``$\res$'' to denote the residue map of $\K$. For any tuple $\bar a=(a_1,\cdots,a_n)$
$\in \K^n$, by $\res(\bar a)$, we mean the tuple $(\res(a_1),\cdots,\res(a_n))$ from the residue field of $\K$. If $X\sq \K^n$, then by $\res(X)$, we mean the set $\{\res(a)|\ a\in X\cap \cO_\K^n\}$.


From now on, We fix $K$ as an elementary submodel of $\K$, namely, $K\equiv \C((t))$. We use $\cO_K$ to denote the valuation ring of $K$, $\kk$ to denote the residue field of $K$, and $\Gamma$ to denote the value group.  $\U(K)$ is the set of units in $\cO_K$. Clearly, $(\cO_K,+)$ and $(\U(K),\times)$ are definable subgroups of the additive group $(K,+)$ and the multiplicative group $(K^*,\times)$, respectively.   If $A\sq B\sq K$ and $p\in S_n(B)$ a complete type over $B$, then by $p|A$ we mean the restriction $\{\varphi\in p |\ \varphi\ \text{is a formula over}\ A\}$. If $p=\tp(\bar a/K)$, then by $\res(p)$ we mean the complete type $\tp(\res(a)/\kk)$.


Our notation for model theory is standard, and we will assume familiarity with basic notions such as very saturated models (or monster models),   partial types, type-definable, finitely satisfiable, etc. References are \cite{Poizat-Book} as well as \cite{Marker-Book}.


\section{Definable subsets of the valuation ring}\label{sec-def-val-ring}
In this section, we will study the definable subsets of $\cO_K$. Let $f(x),g(x)\in K[x]$ be polynomials. Suppose that  
\[
f(x)=b_0+\cdots+b_nx^n.
\]
Let $e_f\in \{b_0,\cdots,b_n\}$  
such that 
\[
v(e_f)=\min\{v(b_i)|\ i=0,\cdots,n\}  
\]
Let $f^*=f/e_f$, then $f^*\in \cO_K[x]$. It is easy to see that $f$ and $f^*$ have the same zeros.  

\begin{lemma}\label{lemm-f(x)=0}
    Suppose that  $a\in \cO_K$ such that $f(a)\in \cO_K$. If $\res(f)(\res(a))\neq 0$, then  $f(a)\neq 0$.
\end{lemma}
\begin{proof}
    We many replace $f$ by $f^*$ since they have the same zeros. Since $\res(f(a))=\res(f)(\res(a))$,  the Lemma follows.
\end{proof}

We conclude directly from the above Lemma that 
\begin{corollary}\label{cor-res(f=0)}
 Let $X=\{a\in \cO_K |\ f(a)=0\}$ and $Z=\{u\in \kk |\ \res(f)(u)=0\}$, then $X\sq \res^{-1}(Z)$ and $\res^{-1}(\kk\backslash Z)\sq \cO_K\backslash X$.
\end{corollary}

\begin{lemma}\label{lemma-v(f)}
    Suppose that  $a\in \cO_K$. If $\res(f^*(a))\neq 0$, then  $v(f(a))=v(e_f)$.
\end{lemma}
\begin{proof}
     Clearly, $v(f(a))=v(e_f)+v(f^*(a))$. Since $\res(f^*(a))\neq 0$, we have $v(f^*(a))=0$. So $v(f(a))=v(e_f)$ as required.
\end{proof}

\begin{corollary}\label{cor-v(f)-leq-v(g)}
  Let $g\in K[x]$.  If $X=\{a\in \cO_K|\ v(f(a))\leq v(g(a))\}$  then there is a finite set $Z\sq \kk$ such that  either $X\sq \res^{-1}(Z)$ or $\res^{-1}(\kk\backslash Z)\sq X$.  
\end{corollary}
\begin{proof}
Suppose that  $g(x)=c_0+\cdots+c_mx^m$. Take  $e_g\in \{c_0,\cdots,c_m\}$ such that $v(e_g)=\min\{v(c_j)|\ j=0,\cdots,m\}$. Let $g^*=g/e_g$ and 
\[
Z=\{c\in \kk|\ (\res(f^*)(c)=0)\vee (\res(g^*)(c)=0)\}.
\]
Then $Z$ is finite. By Lemma \ref{lemma-v(f)}, if $v(e_f)\leq v(e_g)$, then  $\res^{-1}(\kk\backslash Z)\sq X$. If $v(e_f)> v(e_g)$, then  $ X \sq \res^{-1}(Z)$. 
\end{proof}

\begin{lemma}\label{lemma-Pn(f)}
    If $X=\{a\in \cO_K|\ K\models P_n(f(a))\}$, then   there is a finite set $Z\sq \kk$ such that either $X\sq \res^{-1}(Z)$ or $\res^{-1}(\kk\backslash Z)\sq X$.  
\end{lemma}
\begin{proof}
  Let $Z=\{c\in \kk|\ \res(f^*)(c)=0\}$. By Lemma \ref{lemma-v(f)}, $v(f(a))=v(e_f)$ whenever $a\notin \res^{-1}(Z)$. By Corollary \ref{cor-nth-power}, if $v(e_f)\in n\Gamma$, then $\res^{-1}(\kk\backslash Z)\sq X$, otherwise, $X\sq \res^{-1}(Z)$. 
\end{proof}

Let $Z\sq \kk$ and $X\sq \cO_K$, then it is easy to see that $ X \sq \res^{-1}(Z)$ implies $\res^{-1}(\kk\backslash Z)\sq (\cO_K\backslash X)$, and $\res^{-1}(\kk\backslash Z)\sq X$ implies 
$\cO_K\backslash X \sq \res^{-1}(Z)$. 
Summarising  Corollary \ref{cor-res(f=0)},  Corollary \ref{cor-v(f)-leq-v(g)},   and Lemma \ref{lemma-Pn(f)}, we conclude that 
\begin{corollary}\label{coro-res-1(Z)-finite-cofinite}
    Suppose that $X$ is either defined by an $L_{VR}(K)$ formula of \textbf{Type (i)-(iii)}  or negation of a formula of \textbf{Type (i)-(iii)}. Then there is a finite set $Z\sq \kk$ such that either $X\cap \cO_K\sq \res^{-1}(Z)$ or $\res^{-1}(\kk\backslash Z)\sq X\cap \cO_K$. In particular $\res(X\cap \cO_K)$ is either finite or cofinite in $\kk$.
\end{corollary}

\begin{theorem}\label{thm-res(X)-finite-cofinite}
    Let $Y\sq \cO_K$ be definable, there is a finite set $Z\sq \kk$ such that $\res^{-1}(\kk\backslash Z)\sq Y$ if $\res(Y)$ is infinite. In particular, $\res(Y)$ is either finite or cofinite. 
\end{theorem}
\begin{proof}
By quantifier elimination and   Remark \ref{rmk-atomic formulas},  we may assume that $Y$ is a finite boolean combination of sets defined by formulas of \textbf{Type (i)-(iii)}.  Suppose that  $Y=\bigcap_{i=1}^r\bigcup_{j=1}^s Y_{i,j}$, where each $Y_{i,j}$ is either defined by a formula of \textbf{Type (i)-(iii)} or negation of a formula of \textbf{Type (i)-(iii)}.  By Corollary \ref{coro-res-1(Z)-finite-cofinite}, for each $Y_{i,j}$, there is a finite set $Z_{i,j}\sq \kk$ such that either $Y_{i,j}\sq \res^{-1}(Z_{i,j})$ or $\res^{-1}(\kk\backslash Z_{i,j})\sq \cO_K\backslash Y_{i,j}$. Let $Z=\bigcup_{i\leq r,j\leq s}Z_{i,j}$, then we have that either $Y_{i,j}\sq \res^{-1}(Z )$ or $\res^{-1}(\kk\backslash Z )\sq \cO_K\backslash Y_{i,j}$ for each $i\leq r$ and $j\leq s$. Clearly, $\res(Y)\sq \bigcap_{i=1}^r\bigcup_{j=1}^s \res(Y_{i,j})$.


If $\res(Y)$ is infinite, then  for each $i\leq r$ there is $j(i)\leq s$ such that $\res(Y_{i,j(i)})$ is infinite, hence is cofinite in $\kk$. Thus each $\res^{-1}(\kk\backslash Z)$ is contained in $Y_{i,j(i)}$ for each $i\leq r$, we conclude that  so $\res^{-1}(\kk\backslash Z)$ is contained in $Y$. 

\end{proof}

\begin{definition}\label{def-res-finite-cofinite}
    We call a definable subset $X$ of $K$ res-finite (resp. res-cofinite) if $\res(X\cap \cO_K)$ is finite (resp. res-cofinite). We call an $L_{VR}(K)$ formula $\varphi(x)$ is res-finite (resp. res-cofinite) if $\varphi(K)$ is res-finite (resp. res-cofinite).
\end{definition}

By Theorem \ref{thm-res(X)-finite-cofinite}, every definable subset of $K$ is either res-finite or res-cofinite.

\begin{corollary}\label{coro-res-cofinite-meet-every-model}
    Let $X$ be a definable subset of $\cO_K$, $K_0$ an elementary submodel of $K$. If $X$ is res-cofinite, then $X\cap K_0\neq \emptyset$.
\end{corollary}
\begin{proof}
    By Theorem \ref{thm-res(X)-finite-cofinite}, there is a cofinite subset $Z^*\sq \kk$ such that $\res^{-1}(Z^*)\sq X$. Let $\kk_0$ be the residue field of $K_0$. Clearly, $Z^*\cap \kk_0$ is nonempty. Take any $u\in Z^*\cap \kk_0$  and $\Tilde{u}\in \cO_{K_0}$ such that $\res(\Tilde{u})=u$, then $\Tilde{u}\in X\cap K_0$.
\end{proof}


\section{Generic stability of $\cO_{\K}$ and $\U$}\label{sec-Gen-O_k-U}
Recall that by a definable group in $\K$, we mean a definable set $G$ and a definable map $\cdot: G \times G $ such that $(G,\cdot)$ is a group.  We call an $L_{VR}(\K)$-formula $\varphi(\bar x)$  a $G$-formula if $\K\models \forall \bar x(\varphi(\bar x)\to G(\bar x))$. 

Suppose that $G$ is definable over $K$. We can also consider $S_G(K)$ as the space of ultrafilters of the algebra of $G$-formulas over $K$. If $g\in G(K)$ and $\varphi(\bar x)$  is a $G$-formula, then by $g\cdot \varphi$, we mean the formula $\varphi(g^{-1}\cdot \bar x)$. It is easy to see that $(g\cdot \varphi)(K)=g\cdot (\varphi(K))$. If $g\in G(K)$ and $p\in S_G(K)$, then  $g\cdot p=\{g\cdot\varphi|\ \varphi\in p\}$. It is easy to see that $g\cdot p\in S_G(K)$. 

 We call a global type $ p\in S_G(\K)$  generically stable (for $G$) if every $g\cdot p$ is finitely satisfiable in and definable over any elementary submodel of $\K$ for $g\in G$.   We call $G$ generically stable if it admits a global generically stable type $p\in S_G(\K)$. 

In this section, we will prove that the groups $(\cO_{\K},+)$ and $(\U,\times)$ are  generically stable. Let  
    \[
    p_{_{\trans,K}}(x)=\{v(x)=0\}\cup \{h(\res(x))\neq 0| \ h(x)\in \kk[x]\}.
    \]
We see from Remark \ref{res-is-defianble} that $p_{_{\trans,K}}$ is a partial type over $K$.

\begin{lemma}\label{lemma-X-in-p-res(X)-cofinite}
    Let $\varphi(x)$ be an $L_{VR}(K)$-formula, then $  p_{_{\trans,K}}\models \varphi(x)$ iff it is res-cofinite.
\end{lemma}
\begin{proof}
Suppose that $\varphi$ is res-finite. Then 
$\res(\varphi(K)\cap \cO_K)=\{u_1,\cdots,u_n\}$ is finite. Let $h(x)=\prod_{i=1}^n(x-u_i)$, then $h(x)\in \kk[x]$ and  
\[
\{a\in \cO_K|\ h(\res(a))\neq 0\}\cap \varphi(K)=\emptyset.
\]
Namely, $\big(v(x)=0)\wedge (h(\res(x))\neq 0)\big)$ is inconsistent with $\varphi(x)$.
Since 
\[
\big(v(x)=0)\wedge (h(\res(x))\neq 0)\big)\in p_{_{\trans,K}},
\]
we conclude that $ p_{_{\trans,K}}\not\models \varphi(x)$.

 Suppose that $\varphi$ is res-cofinite. Then there is a finite set $Z=\{u_1,\cdots,u_n\}\sq \kk$ such that $\res^{-1}(\kk\backslash Z)\sq \varphi(K)$.  Let $h(x)=\prod_{i=1}^n(x-u_i)$, then 
 \[
 \{a\in K\big(v(a)=0)\wedge (h(\res(a))\neq 0)\}\sq \res^{-1}(\kk\backslash Z)\sq \varphi(K). 
 \]
 so we have
 \[
 K\models \forall x\bigg(\big(v(x)=0)\wedge (h(\res(x))\neq 0)\big)\to\varphi(x) \bigg),
 \]
which implies that  $ p_{_{\trans,K}} \models \varphi(x)$.
\end{proof}

By Theorem \ref{thm-res(X)-finite-cofinite}, an $L_{VR}(K)$-formula $\varphi(x)$ is res-finite (resp. res-cofinite) iff $\ \neg\varphi(x)$ is res-cofinite (resp. res-finite). So we see from  Lemma \ref{lemma-X-in-p-res(X)-cofinite} that $p_{_{\trans, K}}$ determines a complete type over $K$, abusing the notation, this complete type is also denoted by   $p_{_{\trans, K}}$.

The residue field $\kk$ of $K$ is algebraically closed. It is well-known that the theory of algebraically closed fields has quantifier elimination in the language $L_{ring}$. So every definable subset of $\kk$ is finite or cofinite. Let $q_{_{\trans,\kk}}\in S_1(\kk)$ be the unique transcendental type over $\kk$, namely, $q_{_{\trans,\kk}}=\{f(x)\neq 0|\ f\in\kk[x]\}$. It is easy to see that $q_{_{\trans,\kk}}$ is precisely $\res(p_{_{\trans,K}})$: for any $\Tilde{a}\models p_{_{\trans,K}}$,    $\res(\Tilde{a})\models q_{_{\trans,\kk}}$. 

Suppose that $p\in S_n(K)$ is concentrating on $\cO_{\K}^n$ and $q\in S_n(\kk)$. Slightly modifying the definition of ``dominated type''  from \cite{Y.Halevi}, we say that $p$ is dominated by $q$ via the  residue map if for any $a\in \kk^n$,  
\[
a\models q  \iff \Tilde{a}\models p\ \ \text{for all}\  \Tilde{a}\in \res^{-1}(a)
\]
\begin{theorem}\label{thm-p-trans-dominate}
    $p_{_{\trans,K}}$ is dominated by  $q_{_{\trans,\kk}}$ via the residue map. 
\end{theorem}
\begin{proof}
    Suppose that $a\models q_{_{\trans,\kk}}$ and $\Tilde{a}\in \res^{-1}(a)$. Let $\varphi(x)$ be an $L_{VR}(K)$-formula.  By Theorem \ref{thm-res(X)-finite-cofinite}, $\res(\varphi(\K)\cap \cO_{\K})$ is a definable subset of $\kk$. Let $\psi(x)$ be an $L_{ring}$-formula with parameters from $\kk$ such that $\psi(\kk_\K)=\res(\varphi(\K)\cap \cO_{\K})$, where $\kk_\K$ is the residue field of $\K$. Suppose that $\Tilde{a}\models \varphi(x)$, then $a\in \psi(\kk_\K)$. Since $a$ realizes the transcendental type over $\kk$, we see that $\psi(\kk)$ is cofinite. We conclude that  $\Tilde{a}$ realizes every res-cofinite $L_{VR}$-formula over $K$, so $\Tilde{a}$ realizes $p_{_{\trans,K}}$ by Lemma \ref{lemma-X-in-p-res(X)-cofinite}.
\end{proof}

 
We now show that $p_{_{\trans,K}}$ is a generically stable type.

\begin{lemma}\label{lemma-finitely-satisfiable}
 $p_{_{\trans,K}}$ is finitely satisfiable in every  elementary submodel of $K$.
\end{lemma}
\begin{proof}
Let $K_0$ be an elementary submodel of $K$ and $\varphi(x)$ an $L_{VR}(K)$-formula. Suppose that $\varphi(x)\in p_{_{\trans,K}}$.  By Lemma \ref{lemma-X-in-p-res(X)-cofinite},  $\varphi(x)$ is res-cofinite, and by Corollary \ref{coro-res-cofinite-meet-every-model}, $\varphi(K)\cap K_0\neq \emptyset$. This completes the proof.
\end{proof}


\begin{lemma}\label{lemm-p-trans-is definable}
$p_{_{\trans,K}}$ is definable over $\emptyset$.
\end{lemma}
\begin{proof}
   For each $L_{VR}$ formula $\varphi(x,\bar y)$, let  
    \[
    D_{\varphi}=\{\bar b\in K^{|\bar y|}|\ \varphi(x,\bar b)\in p_{_{\trans,K}}\}.
    \] 
    To see the definability of $p_{_{\trans,K}}$, we need to show that: For each $L_{VR}$ formula $\varphi(x,\bar y)$, $D_\varphi$ is $\emptyset$-definable.
    By quantifier elimination and Remark \ref{rmk-atomic formulas}, we only need to check the formulas of \textbf{Type (i)-(iii)}. 
    
    Let $f(x,\bar y )=y_0+y_1x+\cdots+y_nx^n$. By Lemma \ref{lemm-f(x)=0}, we see that   
    \[
    \{\bar b\in K^{n+1}|\ (f(x,\bar b )=0)\in p_{_{\trans,K}}\}=\{(0,\cdots,0)\},
    \]
     which is obviously $\emptyset$-definable. So $D_\varphi$ is $\emptyset$-definable for any formula $\varphi(x,\bar y)$ of \textbf{Type (i)}.

    Suppose that $\varphi(x,y_0,\cdots,y_n,z_0,\cdots,z_m)$ is of the \textbf{Type (ii)}, namely, of the form $v(f(x,\bar y))\leq v(g(x,\bar z))$, where 
    \[
    f(x,\bar y)=y_0+y_1x+\cdots+y_nx^n\]
    and 
    \[
    g(x,\bar z)=z_0+z_1x+\cdots+z_mx^m.
    \]
    By Lemma \ref{lemma-v(f)}, for any $\bar b\in K^{n+1}$ and $\bar c\in K^{m+1}$,  
    \[
    \varphi(x,\bar b,\bar c)\in p_{0,\trans, K}\iff \min\{v(b_i)|\ i\leq n\}\leq \min\{v(c_j)|\ j\leq m\}.
    \]
    Clearly, 
    \[
    D_\varphi=\{(\bar b,\bar c)\in  K^{(n+1)}\times K^{(m+1)}|\ \min\{v(b_i)|\ i\leq n\}\leq \min\{v(c_j)|\ j\leq m\}\}
    \]
    is definable over $\emptyset$. 

    Suppose that $\varphi(x,y_0,\cdots,y_n)$ is of the \textbf{Type (iii)}, namely, of the form $P_n(f(x,\bar y))$, where 
    \[
    f(x,\bar y)=y_0+y_1x+\cdots+y_nx^n
    \]   
    By  Corollary \ref{cor-nth-power} and Lemma \ref{lemma-v(f)}, for each $\bar b=(b_0,\cdots,b_n)\in\K^{n+1}$, we have that   
    \[
    P_n(f(x,\bar b))\in p_{0,\trans, K}\iff  \bigvee_{i=0}^{n}\bigg(\bigwedge_{j=0}^n\big(v(b_i)\leq v(b_j)\big)\to P_n(b_i)\bigg).
    \]
    Clearly, 
    \[
    \{\bar b\in  K^{n+1}|\ \bigvee_{i=0}^{n}\bigg(\bigwedge_{j=0}^n\big(v(b_i)\leq v(b_j)\big)\to P_n(b_i)\bigg)\}
    \]
    is definable over $\emptyset$.
\end{proof}

 Since any definable type over a model has a unique extension, which is also definable, we see from Lemma  \ref{lemm-p-trans-is definable} that:
\begin{corollary}\label{coro-unique-heir}
 For any elementary extension $K_1$ of $K$,   $p_{_{\trans,K_1}}$ is the unique heir of $p_{_{\trans,K}}$ over $K_1$.
\end{corollary}

Now $(\cO_K,+)$ and $(\U(K),\times)$ are definable groups and $p_{_{\trans,K}}\in S_\U(K)\subseteq S_{\cO_\K}(K)$.

\begin{lemma}\label{lemma-p_{0,trans,K}-inv}
For any $a\in \cO_K$, $a+p_{_{\trans,K}}=p_{_{\trans,K}}$. For any $b\in \U(K)$, $b\times p_{_{\trans,K}}=p_{_{\trans,K}}$
\end{lemma}
\begin{proof}
Note that $\res: \cO_K\to \kk$ is a ring homomorphism. Let $A$ be a definable subset of $\cO_K$, $a\in \cO_K$, and $b\in \U(K)$, then $\res(a+A)=\res(a)+\res(A)$ and $\res(b\times A)=\res(b)\times \res(A)$. If $A$ is res-cofinite, then both $a+A$  and $b\times A$   are also res-cofinite. 
By Lemma \ref{lemma-X-in-p-res(X)-cofinite},  we see that 
\[
\varphi(x)\in p_{_{\trans,K}}\iff a+\varphi(x)\in p_{_{\trans,K}}\iff b\times \varphi(x)\in p_{_{\trans,K}}.
\]
This completes the proof.

\end{proof}

Summarizing  Lemma \ref{lemma-finitely-satisfiable}, Lemma \ref{lemm-p-trans-is definable},  and Lemma \ref{lemma-p_{0,trans,K}-inv}, we have that
\begin{theorem}
    $(\cO_{\K},+)$ and $(\U,\times)$ are generically stable, witnessed by a generically stable type $p_{_{\trans,\K}}$.
\end{theorem}

\begin{fact}[See Lemma 9.7 and Lemma 9.12 of \cite{Enrique-Casanovas-book}]
 Let $X\sq \K^n$ be a definable set, $E$ a type-definable equivalence relation on $X$, defined over $K$.  If $|X/E|$ is small/bounded, then for every $a,b\in X$, $a/E=b/E$ whenever $\tp(a/K)=\tp(b/K)$.
\end{fact}

Let $G$ be a group definable in $\K$, $H$ a subgroup group of $G$, and $A\sq \K$ a set of parameters. We call $H$ a type-definable subgroup of $G$ over $A$ if $H\leq G$ is defined by a partial type over $A$. Suppose that $H\leq G$ is a type-definable subgroup over $K$ with small/bounded index, namely, $|G/H|<|\K|$, then the above fact says that each $p\in S_K(G)$ determines a coset of $H$, i.e. , for any $g_1,g_2\in G$, $\tp(g_1/K)=\tp(g_2/K)$ implies that $g_1H=g_2H$. Since $T$ has NIP, by \cite{HPP}, $G$ has the smallest type-definable subgroup of bounded index, written $G^{00}$, which is definable over $\emptyset$ and called the type-definable connected component of $G$. 

\begin{fact}\label{fact-G=G00}
  Let $G$ be a group definable in $\K$.  If there is a global type $p\in S_G(\K)$ which is $G$-invariant, namely, $g\cdot p=p$ for all $g\in G$, then $G=G^{00}$.
\end{fact}
\begin{proof}
    Let $H$ be a type-definable subgroup of $G$ over $A$ of bounded index.  Let $K_0$ be an elementary small submodel of $K$ such that $A\sq K_0$ and $K_0$ meets every coset of $H$. Let $p_0\in S_G(K_0)$ be the restriction of $p$ to $K_0$. Then $g\cdot p_0=p_0$ for all $g\in G(K_0)$. 
    Assume that $p$ is contained in some coset of $H$, say $gH$. If $H$ is a proper subgroup of $G$, then there is $g'\in G$ such that $g'gH\neq gH$. Since $g'\cdot p$ is contained in $g'gH$, we see that $p\neq g'\cdot p$. A contradiction.
\end{proof}

We conclude from Lemma \ref{lemma-p_{0,trans,K}-inv} and Fact \ref{fact-G=G00} that
\begin{corollary}
    Let $G$ be $\cO_{\K}$ or $\U$, then $G=G^{00}$.
\end{corollary}



\section{Generic stability of $\GL(n,\cO_{\K})$}\label{sec-Gen-GLn}
In this section, we identify an $n^2$-tuple $a\in \K^{n\times n}$ with an $n\times n$ matrix.
Recall that $\GL(n,\cO_{\K})=\{g\in \cO_{\K}^{n\times  n}|\ \det(g)\in \U\}$ is the group of invertible $n\times n$ over $\cO_{\K}$. Clearly, $\GL(n,\cO_{\K})$ is definable over $\emptyset$. Suppose that $g=(g_{i,j}) \in \cO_{\K}^{n\times n}$, by $\res(g)$ we mean the matrix $(\res(g_{i,j}))\in \kk_\K^{n\times n}$, where $\kk_\K$ is the residue field of $\K$. Let $\GL(n,\kk_\K)$ be the $n\times n$ general linear algebraic group over $\kk_\K$. It is easy to see that the map:
\[
g\mapsto \res(g),\\ \GL(n,\cO_{\K})\to \GL(n,\kk_\K) 
\]
is an onto group homomorphism. We denote this homomorphism by $\res$ for convenience.
Let $G$ denote the group $\GL(n,\cO_{\K})$ and $\mathfrak G$ denote the group $\GL(n,\kk_\K)$. 

The following fact is easy to verify.
\begin{fact}\label{fact-def-f.s-induction}
Let $M_0\prec M_1\prec M_2\prec M_3$ be structures over a language $L$. Let $b\in M_3$ and $a\in M_2$.
\begin{description}
    \item [(i)] If both $\tp(b/M_2)$ and $\tp(a/M_1)$ are definable over $M_0$, then $\tp(a,b/M_1)$ is definable over $M_0$.
    \item [(ii)] If both $\tp(b/M_2)$ and $\tp(a/M_1)$ are finitely satisfiable in $M_0$, then $\tp(a,b/M_1)$ is finitely satisfiable in $M_0$.
\end{description}
\end{fact}

Now we consider the small submodel $K$. Then 
\[
G(K)=\{g\in \cO_K^{n\times  n}|\ \det(g)\in \U(K)\},
\]
and $\mathfrak{G}(\kk)$ is the $n\times n$ general linear algebraic group over $\kk$. Let 
\[
g^*=(g_1,\cdots g_{n^2})\in G 
\]
and  $K=K_0\prec K_1\prec \cdots \prec K_{n^2}$ an elementary chain such that  $g_i\in K_i$ and  $g_i\models p_{_{\trans, K_{i-1}}}$ for each $i=1,\cdots,n^2$.  Let $\kk_i$ be the residue field of $K_i$ for $i=0,...,n^2$. Since $\res(g_i)$ is transcendental over $\kk_{i}$, for each $i=1,\cdots,n^2$,  the algebraic dimension $\dim(\res(g^*)/\kk)$ of $\res(g^*)$ over $\kk$ is $n^2$, we see that $\res(\det(g^*))=\det(\res(g^*))\neq 0$, so $\det(g^*)\in \U$, and thus $g^*\in G$. Let $p_{_{G,K}}\in S_G(K)$ be the type realized by $g^*$.

\begin{lemma}
  Let $p_{_{G,K}}$ be as the above, then $p_{_{G,K}}$ is definable over and finitely satisfiable in every small submodel of $K$.  
\end{lemma}
\begin{proof}
    Let $M\prec K$. Then by  Lemma \ref{lemma-finitely-satisfiable} and Lemma \ref{lemm-p-trans-is definable},  $\tp(g_i/K_{i-1})$ is definable over and finitely satisfiable in $M$ for each $i=1,\cdots,n^2$. Applying Fact \ref{fact-def-f.s-induction} and induction on $n$, we conclude that $p_{_{G,K}}$ is definable over and finitely satisfiable in $M$.
\end{proof}

    Note that since $p_{_{\trans,K}}$ has a unique heir over any set $A\supseteq  K$, we see that  $p_{_{G,K}}$ is realized by any tuple $h=(h_1,\cdots,h_{n^2})\in G$ such that $h_1\models p_{_{\trans, K}}$, $h_2\models $ the unique heir of $p_{_{\trans, K}}$ over $K\cup\{h_1\}$, $\cdots$, and $h_{n^2}\models$ the unique heir/coheir of $p_{_{\trans, K}}$ over $K\cup\{h_1,\cdots,h_{n^2-1}\}$.

 Note also that there is a unique type $q_{_{\GG,\kk}}\in S_\GG(\kk)$ such that   $a\models q_{_{\GG,\kk}}$ iff $\dim(a/\kk)=n^2$ for each $a\in \GG$ (see Chapter 4 of \cite{Bouscaren-book}  for details), where $\dim(a/\kk)=n^2$ means that $a_1$ is transcendental over $\kk$ and  each $a_i$ is transcendental over $\acl(\kk,a_1,...,a_{i-1})$ for $1<i\leq n$. So $q_{_{\GG,\kk}}=\res(p_{_{G,K}})$, namely, $q_{_{\GG,\kk}}$ is realized by $\res(g^*)$ with $g^*\models p_{_{G,K}}$. We call $q_{{\GG,\kk}}$  the generic type on $\GG$ over $\kk$.  

\begin{theorem}\label{thm-p-G-dominate}
    $p_{_{G,K}}$ is dominated by  $q_{_{\GG,\kk}}$ via the residue map. Namely ,if $a\in \GG$ such that $a\models q_{_{\GG,\kk}}$, then any $\Tilde{a}\in \res^{-1}(a)$ realizes $p_{_{G,K}}$.
\end{theorem}
To prove Theorem \ref{thm-p-G-dominate}, we need to prove the following Lemmas. Let $E$ be any field, then by  $E^{\alg}$, we mean the algebraic closure of $E$. If $E\sq K$, we call that $E$ is algebraically closed in $K$ if $E^{\alg}\cap K=E$.
   
\begin{lemma}\label{lemma-alg-clo-elementary-submodle}
    Suppose that $K_0\sq K$ is algebraically closed in $K$. If there is $u\in K_0$ such that $v(u)=1$,  then $K_0\prec K$.
\end{lemma}
\begin{proof}
Clearly, $K_0$ is a valued subfield of $K$ with valuation ring $\cO_{K_0}=K_0\cap \cO_K$. It is also easy to see that $K_0$ is henselian since $K_0$ is algebraically closed in $K$, and $K$ is henselian. We claim that the residue field $\kk_0$ of $K_0$ is algebraically closed, i.e. an elementary substructure of $\kk$. By Fact \ref{fact-lift}, we can consider $\kk_0\sq \kk$ as  subfields of   $\cO_K$. Since $K_0$ is algebraically closed in $K$, we see that  $\kk_0$ is algebraically closed in $\kk$, which implies that  $\kk_0$ is an algebraically closed field. 

Secondly, we show that $\Gamma_{0}$, the value group  of $K_0$, is an elementary substructure of  $\Gamma$. 
Let $X\sq \Gamma$ be a nonempty set definable over $\Gamma_{0}$. It suffices to show that $X\cap \Gamma_{0}\neq \emptyset$.  By \cite{Raf-Cluckers}, we may assume that $X$ is of the form
\[
X=\{\eta\in \Gamma|\ \alpha<\eta<\beta\wedge D_n(\eta-\gamma) \},
\]
where $\alpha,\beta,\gamma\in \Gamma_{0}$, and $D_n(x)$ is the predicate for ``$x$ is divisible by $n$. Take any $a,b,u\in K_0$ such that $v(a)=\alpha, v(b)=\beta$ and $v(u)=1$. If $b-a \in \N^{>0}$, then there is $i<b-a$ such that $v(u^ia)=i+v(a)\in X\cap \Gamma_0$. Otherwise, assume that  $b-a> \N$. Take any $c\in K_0$ such that $v(c)=\gamma$. There is $i<n$ such that $v(u^ia/c)=i+v(a)-v(c)$ is divisible by $n$. Let $\eta=v(u^ia)$, then $\eta\in X\cap \Gamma_{0}$.

Now we see that $\kk_0\prec \kk$ and $\Gamma_{0}\prec \Gamma$, according to the Ax-Kochen-Ershov-results, $K_0\prec K$ (see Theorem 5.21 of \cite{Dr1} for details).
\end{proof}

\begin{lemma}\label{lemma-lift-ext}
  Let $\kk_\K$ be the residue field of $\K$ and $E$  an algebraically closed subfield of $\kk_\K$. If $E\supseteq \kk$, then there is $K_1\succ K$ such that the residue field of $K_1$ is $E$.
\end{lemma}
\begin{proof}
Let $K_1$ be the algebraic closure of $E\cup K$ in $\K$. By Lemma \ref{lemma-alg-clo-elementary-submodle}, $K_1\prec \K$. Since $K\sq K_1$, we see that $K\prec K_1$. Let $\kk_1$ be the residue field of $K_1$. We now verify that $\kk_1=E$. Clearly, $E\sq \kk_1$. Conversely, consider $\kk_1$ a subfield of $K_1$ and take any $a\in \kk_1$. Let $F$ be the field generated by $E\cup K$. If $a\notin F$, then there is non-constant $f(x)=c_0+\cdots c_nx^n\in F[x]$ such that $f(a)=0$. Let $\lambda=\min\{v(c_i)|\ i=0,\cdots,n\}$ and $f^*=f/\lambda$. Then  $f^*\in E[x]$, $\res(f^*)$ is nonzero, and $\res(f)(a)=\res(f(a))=0$. We conclude that $a$ is algebraic over $E$. Since $E$ is algebraically closed, $a\in E$. This completes the proof.
\end{proof}

\begin{lemma}\label{lemma-lift-sub}
  Let $D$ be an algebraically closed subfield of $\kk$. Then there is $K_0\prec K$ such that $\kk_0=D$, where $\kk_0$ is the residue field of $K_0$.
\end{lemma}
\begin{proof}
    By Fact \ref{fact-lift}, we can consider $D$ as a subfield of $\cO_K$. Take $\alpha\in K$ such that $v(\alpha)=1$. Let $K_0$ be the algebraic closure of $D\cup\{\alpha\}$ in $K$. By Lemma \ref{lemma-alg-clo-elementary-submodle}, $K_0\prec K$. Since $D$ is algebraically closed, a similar argument as in Lemma \ref{lemma-lift-ext} shows that $\kk_0=D$.
\end{proof}


\begin{proof}[Proof of Theorem \ref{thm-p-G-dominate}]
Let $a=(a_1,\cdots,a_{n^2})\models q_{_{\GG,\kk}}$. Let $E_0=\kk$, and each $E_i$ is the algebraic closure of $E_{i-1}(a_{i-1})$ for $i=1,...,n^2$. By Lemma \ref{lemma-lift-ext}, there is an elementary chain $K=K_0\prec K_1\prec...\prec K_{n^2}$ such that $\kk_{i}=E_i$, where $\kk_i$ is the residue field of $K_i$. Since $a_i$ is transcendental over $\kk_{i-1}$,  we have that $a_i\models q_{_{\trans, \kk_{i-1}}}$ for each $i\leq n^2$.  Suppose that $\Tilde{a}=(\Tilde{a}_1,\cdots, \Tilde{a}_{n^2})\in\res^{-1}(a)$. Then by Theorem \ref{thm-p-trans-dominate}, we see that $\Tilde{a}_i\models p_{_{\trans,K_{i-1}}}$ for $i=1,\cdots ,n^2$. So $\Tilde{a}\models p_{_{G,K}}$ as required.
\end{proof}

\begin{corollary}\label{cor-G-inv}
$p_{_{G,K}}$ is $G(K)$-invariant, namely, $g\cdot p_{_{G,K}}=p_{_{G,K}}$ for all $g\in G(K)$. In particular, $p_{_{G,K}}$ is a generically stable type of $G$.
\end{corollary}
\begin{proof}
Let $g\in G(K)$, then $\res(g)\cdot q_{_{\GG,\kk}}=q_{_{\GG,\kk}}$ since $q_{_{\GG,\kk}}$ is the unique generic type of $\GG$ over $\kk$. We see that  if $h\models g\cdot  p_{_{G,K}}$, then $\res(h)\models q_{_{\GG,\kk}}$. Since $p_{_{G,K}}$ is dominated by $q_{_{\GG,\kk}}$ via the residue map,  $g\cdot p_{_{G,K}}=p_{_{G,K}}$ as required.
\end{proof}

\begin{corollary}\label{cor-uniqueness of generically stable}
$p_{_{G,\K}}$ is the unique generically stable type of $G$. 
\end{corollary}
\begin{proof}
Suppose that $r\in S_G(\K)$ is a generically stable type for $G$ and $h\models r$, then $r_0=\tp(\res(h)/\kk_\K)$ is also a generically stable type of $\GG$. Since $\GG$ has a unique generically stable type, we see that $r_0=q_{_{\GG,\kk}}$. Since $p_{_{G,K}}$ is dominated by $q_{_{\GG,\kk}}$ via the residue map, $r=p_{_{G,\K}}$.
\end{proof}

We conclude from Corollary \ref{cor-G-inv} and Corollary \ref{cor-uniqueness of generically stable} that 
\begin{corollary}
    $G=G^{00}$ and $p_{_{G,\K}}$ is the unique generically stable type of $G$.
\end{corollary}

It is easy to see that arguments similar to those above suggest that (linear algebraic) groups over $\cO_\K$ like $\SL(n,\cO_\K)$ and $\SO(n,\cO_\K)$ are generically stable and have unique generically stable types, respectively. Similarly, their generically stable types are dominated by their images via the residue map. It is natural to formulate the following conjecture:
\begin{Conj}
  Let $G\sq \cO_\K^n$ be a $\emptyset$-definable group. Then
  \begin{description}
      \item  [(i)] $G$ is generically stable.
      \item  [(ii)]  $G^{00}$ is a finite index subgroup of $G$.
      \item  [(iii)]  If $p\in S_G(\K)$ is generically stable, then $p$ is dominated by $\res(p)$.      
  \end{description}
\end{Conj} 
\section*{Acknowledgements}
The research is supported by The National Social Science Fund of China (Grant No. 20CZX050).


\begin{thebibliography}{50}
    \bibitem{Y.Halevi}Yatir Halevi. On stably pointed varieties and generically stable groups
in ACVF. Annals of Pure and Applied Logic, 170(2):180–217, 2019.
\bibitem{HP-NIP-inv-measure}Ehud Hrushovski and Anand Pillay. On NIP and invariant measures.
Journal of the European Mathematical Society, 13:1005–1061, 5 2011.
\bibitem{Ax-Kochen}James Ax and Simon Kochen. Diophantine problems over local fields
I. American Journal of Mathematics, 87(3):605–630, 1965.
\bibitem{Ershov}Y. Ershov. On elementary theories of local fields. Algebra i Logika,
pages 5–30, 1965.
\bibitem{Poizat-Book}Bruno Poizat. A Course in Model Theory. Springer New York, NY,
2000.
\bibitem{Marker-Book}David Marker. Model Theory : An Introduction. Springer New York,
NY, 2002.
\bibitem{Enrique-Casanovas-book}Enrique Casanovas. Simple Theories and Hyperimaginaries. Cambridge
University Press, 2011.
\bibitem{HPP}Ehud Hrushovski, Ya’acov Peterzil, and Anand Pillay. Groups, mea-
sures, and the NIP. Journal of the American Mathematical Society,
21(2):563–596, 2008.
\bibitem{Bouscaren-book}Elisabeth Bouscaren. Model Theory and Algebraic Geometry: an introduction to E. Hrushovski's proof of the geometric Mordell-Lang conjecture. Springer Berlin, Heidelberg, 1998.
\bibitem{Raf-Cluckers}Raf Cluckers. Presburger sets and p-minimal fields. The Journal of
Symbolic Logic, 68(1):153–162, 2003.
\bibitem{Dr1}Lou van den Dries. Lectures on the Model Theory of Valued Fields,
pages 55–157. Springer Berlin Heidelberg, Berlin, Heidelberg, 2014.
\bibitem{De1}Françoise Delon. Types sur C((x)). Groupe d’étude de théories stables,
2, 1978-1979.
\end{thebibliography}

\end{document}